\documentclass[runningheads]{llncs}
\usepackage{amssymb}
\usepackage{amsmath}
\usepackage{algorithm}
\usepackage{algpseudocode}
\usepackage{cite}
\usepackage{mathtools}
\DeclareMathOperator*{\minimize}{minimize}
\usepackage{multirow}

\usepackage{graphicx}
\DeclareGraphicsExtensions{.pdf,.png}
\begin{document}
\title{Atomic Super-Resolution Tomography}
%
%
\author{Poulami Somanya Ganguly\inst{1,2}\orcidID{0000-0003-1672-7194} \and
Felix Lucka\inst{1,3}\orcidID{0000-0002-8763-5177} \and
Hermen Jan Hupkes\inst{2}\orcidID{0000-0003-1726-5323} \and
Kees Joost Batenburg\inst{1,2}\orcidID{0000-0002-3154-8576}
}
\authorrunning{P.S.~Ganguly, F.~Lucka, H.J.~Hupkes, K.J.~Batenburg}
%
\institute{Centrum Wiskunde \& Informatica, Science Park, Amsterdam, The Netherlands
\and
The Mathematical Institute, Leiden University, Leiden, The Netherlands
\and
Centre for Medical Image Computing, University College London, London, United Kingdom
\\
\email{\{poulami.ganguly, felix.lucka, joost.batenburg\}@cwi.nl}
\email{hhupkes@math.leidenuniv.nl}}
%
\maketitle              
\begin{abstract}
We consider the problem of reconstructing a nanocrystal at atomic resolution from electron microscopy images taken at a few tilt angles. A popular reconstruction approach called discrete tomography confines the atom locations to a coarse spatial grid, which is inspired by the physical \emph{a priori} knowledge that atoms in a crystalline solid tend to form regular lattices. Although this constraint has proven to be powerful for solving this very under-determined inverse problem in many cases, its key limitation is that, in practice, defects may occur that cause atoms to deviate from regular lattice positions. Here we propose a grid-free discrete tomography algorithm that allows for continuous deviations of the atom locations similar to super-resolution approaches for microscopy. The new formulation allows us to define atomic interaction potentials explicitly, which results in a both meaningful and powerful incorporation of the available physical \emph{a priori} knowledge about the crystal's properties. In computational experiments, we compare the proposed grid-free method to established grid-based approaches and show that our approach can indeed recover the atom positions more accurately for common lattice defects.

\end{abstract}
\section{Introduction}
 Electron tomography is a powerful technique for resolving the interior of nanomaterials. After preparing a microscopic sample, a series of projection images (so called tilt-series) is acquired by rotating the specimen in the electron microscope, acquiring data from a range of angles. In recent years, electron tomography has been successfully applied to reconstruct the 3D positions of the individual atoms in nanocrystalline materials \cite{van2012big, chen2013three, rez2013three}.

Since the first demonstration of atomic resolution tomography of nanocrystals in 2010 by discrete tomography \cite{van2011three}, a range of tomographic acquisition techniques and reconstruction algorithms have been applied to reconstruct nanocrystals of increasing complexity. In the discrete tomography approach, atoms are assumed to lie on a regular lattice and the measured projections can be considered as atom counts along lattice lines. A key advantage of this approach is its ability to exploit the constraints induced by the discrete domain and range of the image. As a consequence, a small number of projection angles (typically less than 5) can already lead to an accurate reconstruction \cite{batenburg2006,batenburg2009}. The theoretical properties of the discrete reconstruction problem have been studied extensively with results on algorithm complexity, uniqueness, and stability \cite{gardner1997,baake2006,alpers2006}. A key drawback of the discrete lattice assumption when considering real-world applications to nanocrystal data is that in many interesting cases the atoms do not lie on a perfect lattice due to defects in the crystal structure or interfaces between different crystal lattices. In such cases the atoms do not project perfectly into columns, forming a mismatch with the discrete tomography model.

As an alternative, it has been demonstrated that a more conventional tomographic series consisting of hundreds of projections of a nanocrystal can be acquired in certain cases. An image of the nanocrystal is then reconstructed using sparsity based reconstruction techniques on a continuous model of the tomography problem. This approach does not depend on the lattice structure and allows one to reconstruct defects and interfaces \cite{goris2015measuring}. As a downside, the number of required projections is large and to accurately model the atom positions the reconstruction must be represented on a high-resolution pixel grid resulting in a large-scale computational problem. This raises the question if a reconstruction problem can be defined that fills the gap between these two extremes and can exploit the discrete nature of the lattice structure while at the same time allowing for continuous deviations of atom positions from the perfect lattice.

In this paper we propose a model for the atomic resolution tomography problem that combines these two characteristics. Inspired by the algorithm proposed in \cite{boyd2017alternating}, our model is based on representing the crystal image as a superposition of delta functions with continuous coordinates and exploiting sparsity of the image to reduce the number of required projections. We show that by incorporating a physical model for the potential energy of the atomic configuration, the reconstruction results can be further improved. 

\section{Problem Setting} \label{mathsetup}
In this section we formulate a mathematical model of the atomic resolution tomography problem and discuss several approaches to solve it. Some of these approaches assume that the atom locations are restricted to a perfect grid, the \textit{crystal lattice}, which corresponds to only one possible local minimum of the potential energy of the atomic configuration. To overcome certain limitations of this assumption, we propose an alternative formulation where the atom locations are allowed to vary continuously and an explicit model of the potential energy of their configuration is used to regularize the image reconstruction. 

An atomic configuration is characterized by a positive measure $\mu$ on a bounded subset $X$ of $\mathbb{R}^d$. We denote the space of such measures by $\mathcal{M}(X)$. The measure represents the \emph{electron density}, which is the probability that an electron is present at a given location. The electron density around an atomic configuration is highest in regions where atoms are present. In electron tomography, electron density is probed by irradiating a sample with a beam of electrons. The beam undergoes absorption and scattering due to its interactions with the electrons of the atomic configuration. The transmitted or scattered signal can then be used to form an image. The Radon transform provides a simplified mathematical model of this ray-based image formation process. For $d=2$, the Radon transform $\mathcal{R} \mu$ can be expressed as integrals taken over straight rays
\begin{align} \label{eq:radon}
\mathcal{R}[\mu](r, \theta) &:= \int_{l(r,\theta)} d\mu,  \\
l(r,\theta) &= \{(x_1,x_2) \in \mathbb{R}^2 \,|\, x_1\cos\theta + x_2\sin\theta = r \}, 
\end{align}
where we parametrized the rays by the projection angle $\theta$ and the distance on the detector $r$. The corresponding inverse problem is to recover $\mu$ from noisy observations of $y = \mathcal{R} \mu + \varepsilon$. One way to formulate a solution to this problem is via the following optimisation over the space of measures:
\begin{equation}\label{eq:opt0}
    \minimize_{\mu \in \mathcal{M}(X)} \quad \| \mathcal{R} \mu - y \|^2_2 , 
\end{equation}
which is an infinite dimensional non-negative linear least-squares problem. In the following, we will introduce a series of discretisations of this optimisation problem. Numerical schemes to solve them will be discussed in Section \ref{algos}. 

In situations where we only have access to data from a few projection angles, introducing a suitable discretisation of \eqref{eq:opt0} is essential for obtaining a stable reconstruction. One way to achieve this is to restrict the atom locations to a spatial grid with $n$ nodes, ${\vec{x}}_{i=1}^n$, and model their interaction zone with the electron beam by a Gaussian with known shape $G$. The atom centres are then delta peaks $\delta_{\vec{x}_i}$ on the gridded image domain. The Gaussian convolution of atom centres can be viewed as the ``blurring" produced by thermal motion of atoms. In fact, it is known from lattice vibration theory that, for large configurations, the probability density function of an atom around its equilibrium position is a Gaussian, whose width depends on temperature, dimensionality and interatomic forces \cite{montroll1956theory}. The discretized measure $\mu$ can then be written as
\begin{equation}\label{eq:gridmu}
    \mu_{\text{grid}} = \sum_{i = 1}^n w_i (G \ast \delta_{\vec{x}_i}),
\end{equation}
where $n$ is the total number of grid points and weights $w_i \geq 0$ were introduced to indicate confidence in the presence or absence of an atom at grid location $i$. If we insert \eqref{eq:gridmu} in \eqref{eq:opt0} and introduce the forward projection of a single atom as $\psi_i \coloneqq \mathcal{R}(G \ast \delta_{\vec{x}_i})$ we get
\begin{multline*}
   \| \mathcal{R} \mu_{\text{grid}} - y \|^2_2 =  \| \mathcal{R} \sum_{i = 1}^n w_i (G \ast \delta_{\vec{x}_i}) - y \|^2_2 =  \\
    \| \sum_{i = 1}^n w_i \mathcal{R} (G \ast \delta_{\vec{x}_i}) - y \|^2_2 =:  \| \sum_{i = 1}^n \psi_i w_i  - y \|^2_2 =:  \| \mathrm{\Psi} w  - y \|^2_2 
\end{multline*}
The corresponding optimisation problem is given by 
\begin{equation} \label{eq:opt1}
 \minimize_{w \in \mathbb{R}^n_+} \quad \| \mathrm{\Psi} w  - y \|^2_2 \quad , 
\end{equation}
which is a finite dimensional linear non-negative least squares problem. 

The choice of the computational grid in \eqref{eq:gridmu} is unfortunately not trivial. Only in certain situations, one can assume that all atoms lie on a lattice of known grid size and orientation, and directly match this lattice with the computational grid. In general, one needs to pick a computational grid of much smaller grid size. With the data $y$ given, the grid admits multiple solutions of \eqref{eq:opt1} and most efficient computational schemes tend to pick a blurred, artefact-ridden solution with many non-zero weights far from the true, underlying $\mu$, as we will demonstrate in Section \ref{expresults}. To obtain a better reconstruction, one can choose to add \emph{sparsity constraints} which embed our physical \emph{a priori} knowledge that $\mu$ originates from a discrete configuration of atoms. In our model \eqref{eq:gridmu}, this corresponds to a $w \in \mathbb{R}^n_+$ with few non-zeros entries. To obtain such a sparse solution we can add a constraint on $\ell^0$ norm of the weights to the optimisation problem:
\begin{equation*}
    \begin{aligned}
    & \underset{w \in \mathbb{R}^n_{+}}{\minimize}
    & & \| \mathrm{\Psi}w - y \|^2_2 \\
    & \text{subject to}
    & & |w|_0 \leq K.
    \end{aligned}
\end{equation*}
However, this problem is NP-hard \cite{foucart2013mathematical}. A approximate solution can be found by replacing the $\ell^0$ norm with the $\ell^1$ norm and adding it to the objective function:
\begin{equation}\label{eq:opt2}
        \minimize_{w \in \mathbb{R}^n_{+}} \quad \| \mathrm{\Psi}w - y \|^2_2 + \lambda \|w\|_1,
\end{equation}
where $\lambda$ is the relative weight of the sparsity-inducing term. This particular choice of formulation is not always best and alternative formulations of the same problem exist \cite{foucart2013mathematical}.

For atomic configurations where only one type of atom is present, the weights can be considered to be one where an atom is present and zero everywhere else. This corresponds to discretising the range of the reconstructed image. The fully discrete optimisation problem then becomes:
\begin{equation}\label{eq:opt3}
\minimize_{w \in \{0,1\}^n} \quad \| \mathrm{\Psi}w - y \|^2_2.
\end{equation}
With image range discretisation, a constraint on the number of atoms is typically no longer needed because adding an additional atom with weight 1 after all atoms have been found leads to an increase in the objective function. 

\begin{figure}
    \centering
    \includegraphics[width=\textwidth]{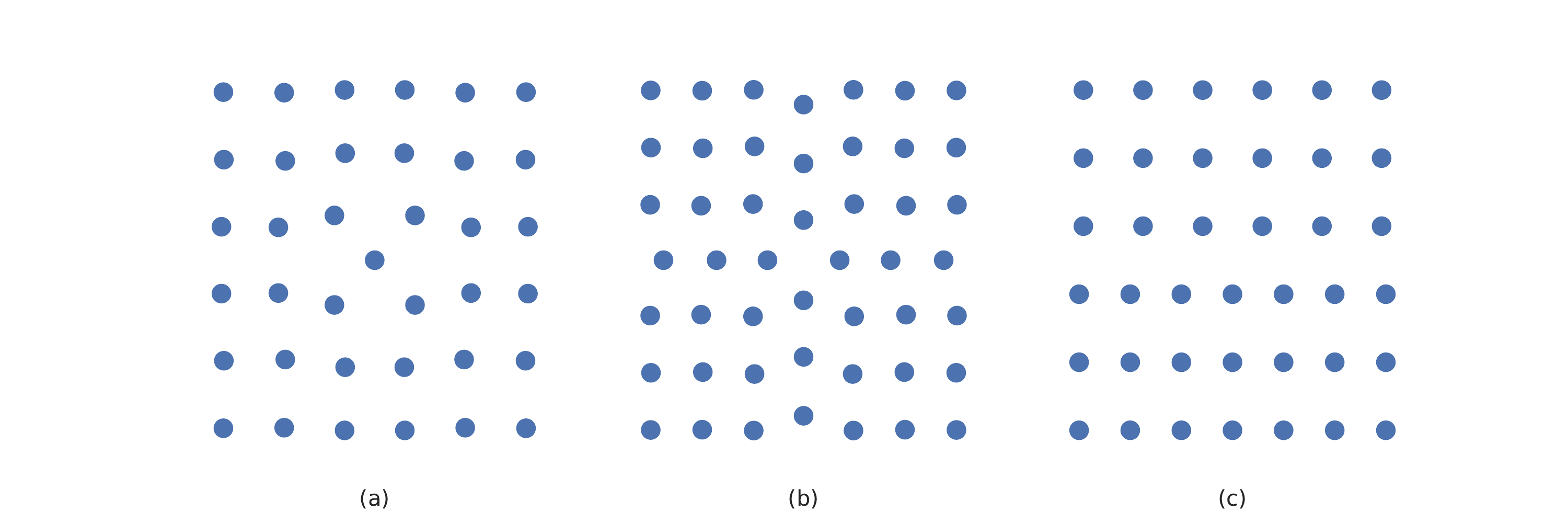}
    \caption{Atomic configurations with (a) an interstitial point defect, (b) a vacancy and (c) an edge dislocation.}
    \label{fig:configs}
\end{figure}

Although the optimisation problems (\ref{eq:opt1}), (\ref{eq:opt2}) and (\ref{eq:opt3}) allow for the recovery of atomic configurations without solving (\ref{eq:opt0}), all of them rely crucially on discretisation of the domain of the reconstructed image, i.e.~the assumption that atoms lie on a grid. However, this assumption is not always true. In particular, atomic configurations often contain defects where atom positions deviate from the perfect lattice. Fig.~\ref{fig:configs} shows examples of common lattice defects. In order to resolve these defects correctly, the image domain must be discretized to higher resolutions, i.e.~the grid of possible atom positions must be made finer. This introduces two main problems: First, making the grid finer for the same data makes the inverse problem more ill-posed. Second, the computational time increases significantly even for modestly sized configurations.

In order to overcome these difficulties, we revisit (\ref{eq:gridmu}) and remove the requirement for $\vec{x}_i$ to lie on a grid. The projection of a single atom now becomes a function of its location $\vec{x} \in \mathbb{R}^d$, $\psi(\vec{x}) := \mathcal{R}(G \ast \delta_{\vec{x}})$. We keep the image range discretisation introduced above by requiring $w_i \in \{0,1\}$. Now, \eqref{eq:opt3} becomes 
\begin{equation}\label{eq:infoptnonconv}
    \begin{aligned}
    & \underset{\vec{x}\in X^n, w \in \{0,1\}^n}{\minimize}
    & & \Big \| \sum_{i=1}^n w_i \psi(\vec{x}_i) - y \Big \|^2_2. \\
    \end{aligned}
\end{equation}
The minimisation over $\vec{x}$ is a non-linear, non-convex least-squares problem which has been studied extensively in the context of mathematical super-resolution \cite{candes2014towards, bredies2013inverse, boyd2017alternating}. In these works, efficient algorithms are derived from relating it back to the infinite dimensional linear least-squares problem on the space of measures \eqref{eq:opt0}. For instance, for applications such as fluorescence microscopy \cite{boyd2017alternating} and ultrasound imaging \cite{alberti2019dynamic}, an alternating descent conditional gradient (ADCG) algorithm has been proposed, which we will revisit in the next section. Compared to these works, we have a more complicated non-local and under-determined inverse problem and the minimisation over $w$ adds a combinatorial, discrete flavor to \eqref{eq:infoptnonconv}. To further tailor it to our specific application, we will incorporate physical \emph{a priori} knowledge about atomic configurations of crystalline solids by adding a functional formed by the atomic interaction potentials. This will act as a regularisation of the underlying under-determined inverse problem.

\begin{figure}[t]
    \centering
    \includegraphics[width=0.8\textwidth]{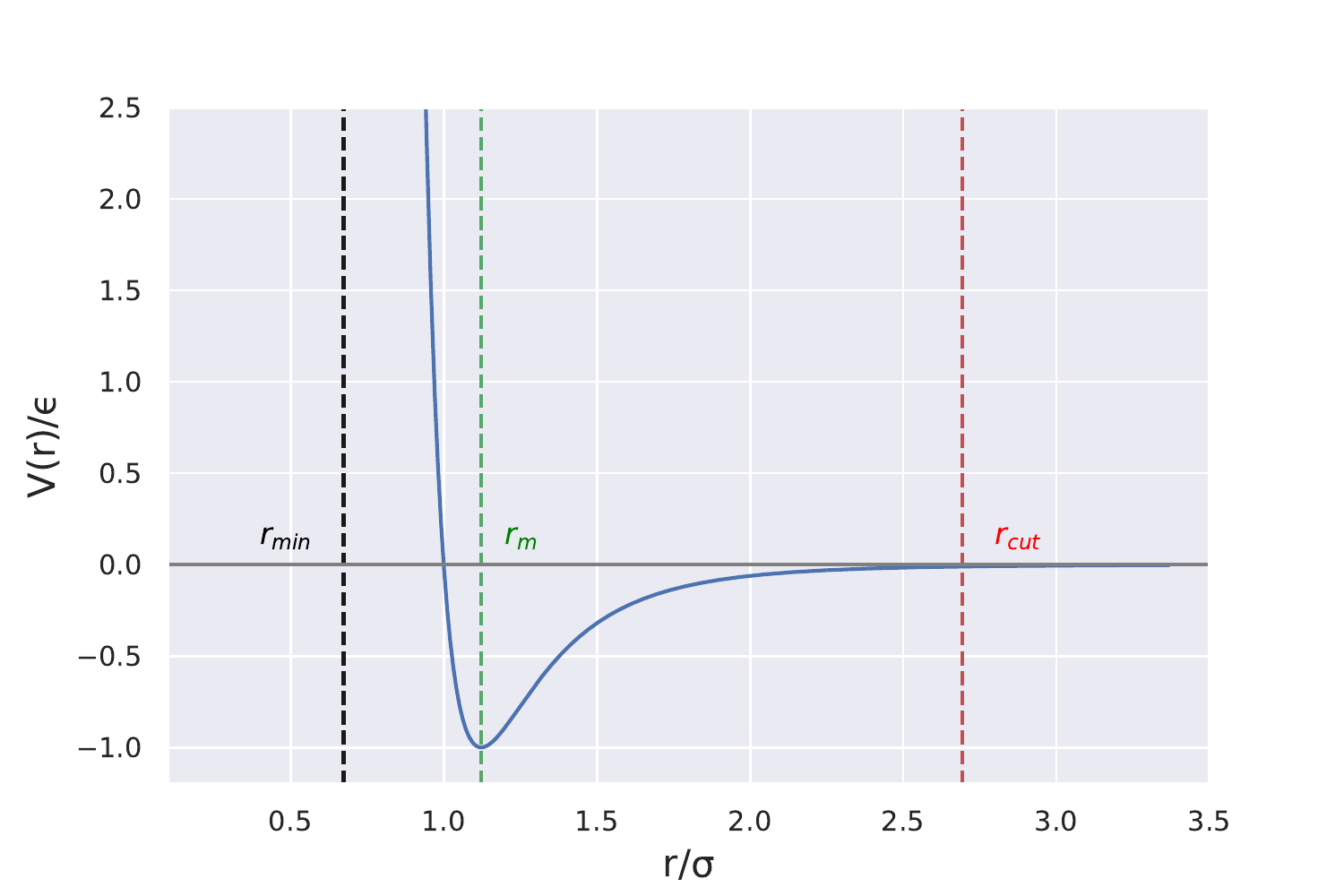}
    \caption{The normalized Lennard-Jones pair potential as a function of normalized interatomic separation.}
    \label{fig:LJ}
\end{figure}

\subsection{Potential Energy of the Atomic Configuration}
The total energy of an atomic configuration is the sum of its potential energy and kinetic energy. As we consider only static configurations, the kinetic energy of the configuration is zero and the total energy is equal to the potential energy. In order to compute the potential energy of the atomic configuration, we must prescribe the interaction between atoms. In this paper, we use the Lennard-Jones pair potential, which is a simplified model of interatomic interactions. The Lennard-Jones potential $V_{\text{LJ}}$ as a function of interatomic separation $r$ is given by \cite{frenkel2001understanding}
\begin{align}\label{eq:lj}
    V_{\text{LJ}}(r) &= 
    \begin{cases}4 \epsilon \Big[ \Big(\dfrac{\sigma}{r}\Big)^{12} - \Big(\dfrac{\sigma}{r}\Big)^6 \Big], & r<r_{\mathrm{cut}}\\
    0, & r \geq r_{\mathrm{cut}}
    \end{cases}
\end{align}
where $\epsilon$ is the depth of the potential well and $\sigma$ is the interatomic separation at which the potential is zero. The separation at which the potential reaches its minimum is given by $r_m = 2^{1/6}\sigma$. The parameter $r_{\mathrm{cut}}$ denotes a cut-off separation beyond which the potential is inactive. Fig.~\ref{fig:LJ} shows the form of the the Lennard-Jones pair potential as a function of interatomic separation. The potential energy of the atomic configuration is computed by summing over the pairwise interaction between all pairs of atoms 
\begin{equation}
V_{\text{tot}}(\vec{x}_1, \vec{x}_2,...,\vec{x}_N) = 
\sum_{i>j} V_{\text{LJ}}(\vec{x}_i - \vec{x}_j).
\end{equation}
Adding this energy to the objective in (\ref{eq:infoptnonconv}) leads to 
\begin{equation}\label{eq:infoptwithenergy}
    \minimize_{\vec{x}\in \mathcal{C}, w \in \{0,1\}^n} \quad \Big \| \sum_{i=1}^n w_i \psi(\vec{x}_i) - y \Big \|^2_2 + \alpha V_{\text{tot}}(\vec{x}).
\end{equation}
The regularisation parameter, $\alpha$, adjusts the relative weight of the energy term, so that by tuning it we are able to move between atomic configurations that are data-optimal and those that are energy-optimal. The constraint set $\mathcal{C} \subset X^n$ is defined by a minimum distance $r_{\mathrm{min}}$, such that $|\vec{x}_i - \vec{x}_j|>r_{\mathrm{min}}$, $\forall i>j$. The minimum distance, $r_{\mathrm{min}}$, is chosen to be smaller than the optimal interatomic separation $r_m$ and allows us to set $\alpha$ to 0 and still avoid configurations where atoms are placed exactly at the same position. For small separations, the energy is dominated by the $\big(\frac{\sigma}{r}\big)^{12}$ term and increases sharply for separations less than $r_m$. Thus, for non-zero $\alpha$, configurations where atoms have a separation less than $r_m$ are highly unlikely.

\section{Algorithms}\label{algos}
In this section we discuss several algorithms to solve the optimisation problems introduced in Section \ref{mathsetup}.

\subsection{Projected Gradient Descent}\label{sirt} The non-negative least-squares problem (\ref{eq:opt1}) can be solved with a simple iterative first-order optimisation scheme. At each step of the algorithm, the next iterate is computed by moving in the direction of the negative gradient of the objective function. Non-negativity of the weights is enforced by projecting negative iterates to zero. Mathematically, each iterate is given by
\begin{align} \label{projGrad}
    w^{k+1} = \sideset{}{_+}\prod \Big (w^k + t \mathrm{\Psi}^T (\mathrm{\Psi}w^k-y)\Big),
\end{align}
where $t$ is the step size and the projection operator is given by
\begin{equation} \label{proj}
    \sideset{}{_+}\prod(\cdot) = \text{max}(\,\cdot\, ,0).
\end{equation}
In the numerical experiments in Section \ref{expresults}, we used the SIRT algorithm \cite{kak2002principles} as implemented in the tomographic reconstruction library ASTRA \cite{van2015astra}, which is based on a minor modification of the iteration described above.  

\subsection{Proximal Gradient Descent}\label{prox}
If we add the non-smooth $\ell^1$ regularizer and obtain problem (\ref{eq:opt2}), we need to extend \eqref{projGrad} to a proximal gradient scheme \cite{parikh2014proximal} 
\begin{equation} \label{proxGrad}
    w^{k+1} = \text{prox}_{h}\Big (w^k + t \mathrm{\Psi}^T (\mathrm{\Psi}w^k-y)\Big), 
\end{equation}
where the projection operator \eqref{proj} is replaced by the proximal operator of the convex function 
\begin{equation}
    h(x) := \begin{cases}
            \lambda \|x\|_1  & x \geq 0 \\
            0                & \text{elsewhere}
    \end{cases} \quad ,
\end{equation}
which is given by the non-negative soft-thresholding operator
\begin{align*}
    \text{prox}_{h}(x) &= 
    \begin{cases}
    x - \lambda, & x \geq \lambda \\
    0, & \text{elsewhere}
    \end{cases} \quad .
\end{align*}
In the numerical experiments in Section \ref{expresults}, we used the fast iterative soft-thresholding algorithm (FISTA) \cite{beck2009fast} as implemented in the Python library ODL \cite{adler2017operator}, which is based on a slight modification of the iteration described above.  


	

\begin{algorithm}
	\caption{Discrete simulated annealing} 
	\begin{algorithmic}[1]
	    \While {$\beta < \beta_{\text{max}}$}
		\State Select new atom location: $\tilde{w}^{k} \in \arg\min_{k \in \mathcal{C}} \, \mathrm{\Psi} w^k - y$
		\State Add new atom to current configuration: $\tilde{w}^{k+1} \leftarrow \{w^k,\tilde{w}^k\}$
		\State Accept new configuration with a certain probability:
		\If {$\beta\|\mathrm{\Psi}\tilde{w}^{k+1}-y\|^2_2 < \beta\|\mathrm{\Psi}w^{k}-y\|^2_2$} 
		\State $w^{k+1} \leftarrow \tilde{w}^{k+1}$
		\Else
		\State Generate random number: $t \in \text{rand}[0,1)$
		\If {$t < e^{-\beta\|\mathrm{\Psi}\tilde{w}^{k+1}-y\|^2_2}/e^{-\beta\|\mathrm{\Psi}{w}^{k}-y\|^2_2}$}
		\State $w^{k+1} \leftarrow \tilde{w}^{k+1}$
		\EndIf
		\EndIf
		\State Move atom: $w^{k+1} \leftarrow \text{random move}(w^{k+1})$
		\State Run acceptance steps 5--13
		\State Increase $\beta$
	\EndWhile
	\end{algorithmic} 
	\label{alg:sa}
\end{algorithm}

\subsection{Simulated Annealing}\label{sa} For solving the fully discrete problem (\ref{eq:opt3}), we used a simulated annealing algorithm as shown in Algorithm \ref{alg:sa}, which consists of two subsequent accept-reject steps carried out with respect to the same inverse temperature parameter $\beta$. In the first one, the algorithm tries to add a new atom to the existing configuration. In the second one, the atom locations are perturbed locally. As $\beta$ is increased towards $\beta_{\text{max}}$, fewer new configurations are accepted and the algorithm converges to a minimum. 
  
In the atom adding step at each iteration $k$, the algorithm tries to add an atom at one of the grid location $i$ where the residual $\mathrm{\Psi} w^k - y$ is minimal (this corresponds to flipping $w^k_i$ from $0$ to $1$ in (\ref{eq:opt3})). The allowed grid locations belong to a constraint set $\mathcal{C}$, such that no two atoms are closer than a pre-specified minimum distance $r_{\mathrm{min}}$. 
To perturb the atom positions locally, the algorithm selects an atom at random and moves it to one of its 4 nearest neighbor locations at random. 

\subsection{ADCG with Energy}\label{adcgenergy}

\begin{algorithm}[b]
	\caption{ADCG with energy} 
	\begin{algorithmic}[1]
	\For {$k=1:k_{\text{max}}$}
            \State Compute next atom in grid $g$: 
            \Statex \hskip \algorithmicindent $\vec{x}_{\mathrm{new}} \in \arg \min_{\vec{x}_{\mathrm{new}} \in g, (\vec{x}^k,\vec{x}_{\mathrm{new}}) \in \mathcal{C}} \| \sum_{i=1}^k \psi(\vec{x}_i) - y + \psi(\vec{x}_{\mathrm{new}})\| + \alpha V_{\text{tot}}(\vec{x}^k,\vec{x}_{\mathrm{new}})$
            \State Update support: $\vec{x}^{k+1} \leftarrow \{ \vec{x}^k, \vec{x}_{\mathrm{new}} \}$
            \State Locally move atoms:
	\Statex \hskip \algorithmicindent $\vec{x}^{k+1} \leftarrow \min_{\vec{x} \in X} \| \mathrm{\Psi} \mu(\vec{x}^{k+1}) - y\|^2_2 + \alpha V(\vec{x}^{k+1})$
	\State Break if objective function is increasing:
	\If {$\| \mathrm{\Psi} \mu(\vec{x}^{k+1}) - y\|^2_2 + \alpha V(\vec{x}^{k+1}) > \| \mathrm{\Psi} \mu(\vec{x}^{k}) - y\|^2_2 + \alpha V(\vec{x}^{k})$} break
	\EndIf
	\EndFor
	\end{algorithmic} 
	\label{alg:adcg_energy}
\end{algorithm}

Variants of the Frank-Wolfe algorithm (or conditional gradient method) \cite{frank1956algorithm,jaggi2013revisiting} have been proposed for solving problems of the form \eqref{eq:infoptnonconv} \cite{denoyelle2019sliding, alberti2019dynamic} without discrete constraints for $w$ and are commonly known as alternating descent conditional gradient (ADCG) schemes (see \cite{poon2019multidimensional} for an analysis specific to multidimensional sparse inverse problems). Here, we modify the ADCG scheme to 
\begin{enumerate}
    \item incorporate binary constraints on $w$
    \item handle the singularities of the atomic interaction potentials 
    \item avoid local minima resulting from poor initialisations
\end{enumerate}{}
The complete algorithm is shown in Algorithm \ref{alg:adcg_energy}. Essentially, the scheme also alternates between adding a new atom to the current configuration and optimising the positions of the atoms. 

In the first step, the image domain is coarsely gridded and the objective function after adding an atom at each location is computed. Locations closer to existing atoms than $r_{\mathrm{min}}$ are excluded. In the second step, the atom coordinates are optimized by a continuous local optimisation method. Here, the Nelder-Mead method \cite{nelder1965simplex} implemented in SciPy \cite{virtanen2019scipy} was used. 

A continuation strategy is used to avoid problems resuling from poor initilisations: Algorithm \ref{alg:adcg_energy} is run for increasing values of $\alpha$, starting from $\alpha=0$. The reconstruction obtained at the end of a run is used as initialisation for the next. In the following section, we demonstrate the effect of increasing $\alpha$ on the reconstructions obtained and discuss how an optimal $\alpha$ was selected. In the following section, we refer to  Algorithm \ref{alg:adcg_energy} as ``ADCG" when used for $\alpha = 0$ and as ``ADCG with energy" otherwise.

\section{Numerical Experiments}\label{expresults}
We conducted numerical experiments by creating 2D atomic configurations with defects and using the algorithms discussed in Section \ref{algos} to resolve atom positions. In this section we describe how the ground truth configurations were generated and projected, and compare the reconstruction results of different algorithms.

\subsection{Ground Truth Configurations}
We generated ground truth configurations using the molecular dynamics software HOOMD-blue \cite{anderson2008general, glaser2015strong}. We created perfect square lattices and then induced defects by adding or removing atoms. The resulting configuration was then relaxed to an energy minimum using the FIRE energy minimizer \cite{bitzek2006structural} to give the configurations shown in Fig.~\ref{fig:configs}. The following parameter values were used in \eqref{eq:lj} for specifying the Lennard-Jones pair potential between atoms.

\begin{center}
 \begin{tabular}{ |p{3cm}||p{3cm}|p{3cm}|p{3cm}|  }
 \hline
 Defect type & $\epsilon$ & $\sigma$ & $r_{\mathrm{cut}}$ \\ [1.5ex] 
 \hline\hline
 Interstitial defect & 0.4 & 0.15 & 0.4 \\ 
 \hline
 Vacancy & 0.4 & 0.14 & 0.4 \\
 \hline
 Edge dislocation & 0.4 & 0.13 & 0.17 \\
 \hline
\end{tabular}
\end{center}

\subsection{Discretized Projection Data}
We generated two 1D projections for each ground truth atomic configuration at projection angles $\theta = 0^\circ, 90^\circ$. As discussed in Section \ref{mathsetup}, the projection of a single atom centre is given by a Gaussian convolution followed by the Radon transform. The Radon transform of a Gaussian is also a Gaussian. Therefore, we interchanged the two operations in the forward transform to speed up the computations. The sum over individual projections of atom centres was used as the total (noise-free) projection. Using the Radon transform in \eqref{eq:radon}, each atom centre was projected onto a 1D detector, following which it was convolved with a 1D Gaussian of the form $G(z) = e^{-(z - z_0)/\varsigma^2}$, where $z_0$ is the position of the atom centre on the detector and $\varsigma$ controls the width of the Gaussian. Finally, the continuous projection was sampled at a fixed number of points to give rise to a \emph{discrete} projection. For our experiments, the $\varsigma$ of the Gaussian function was taken to be equal to the discretisation of the detector given by the detector pixel size $d$. Both were taken to be 0.01.

\subsection{Discretisation of Reconstruction Volume}
For SIRT, FISTA and simulated annealing (described in subsections \ref{sirt}, \ref{prox} and \ref{sa}, respectively), each dimension of the reconstruction area was discretized using the detector pixel size $d$. Therefore, there were $1/d \times 1/d$ grid points in total. 


\begin{figure}[H]
  \centering
    \includegraphics[scale=0.24]{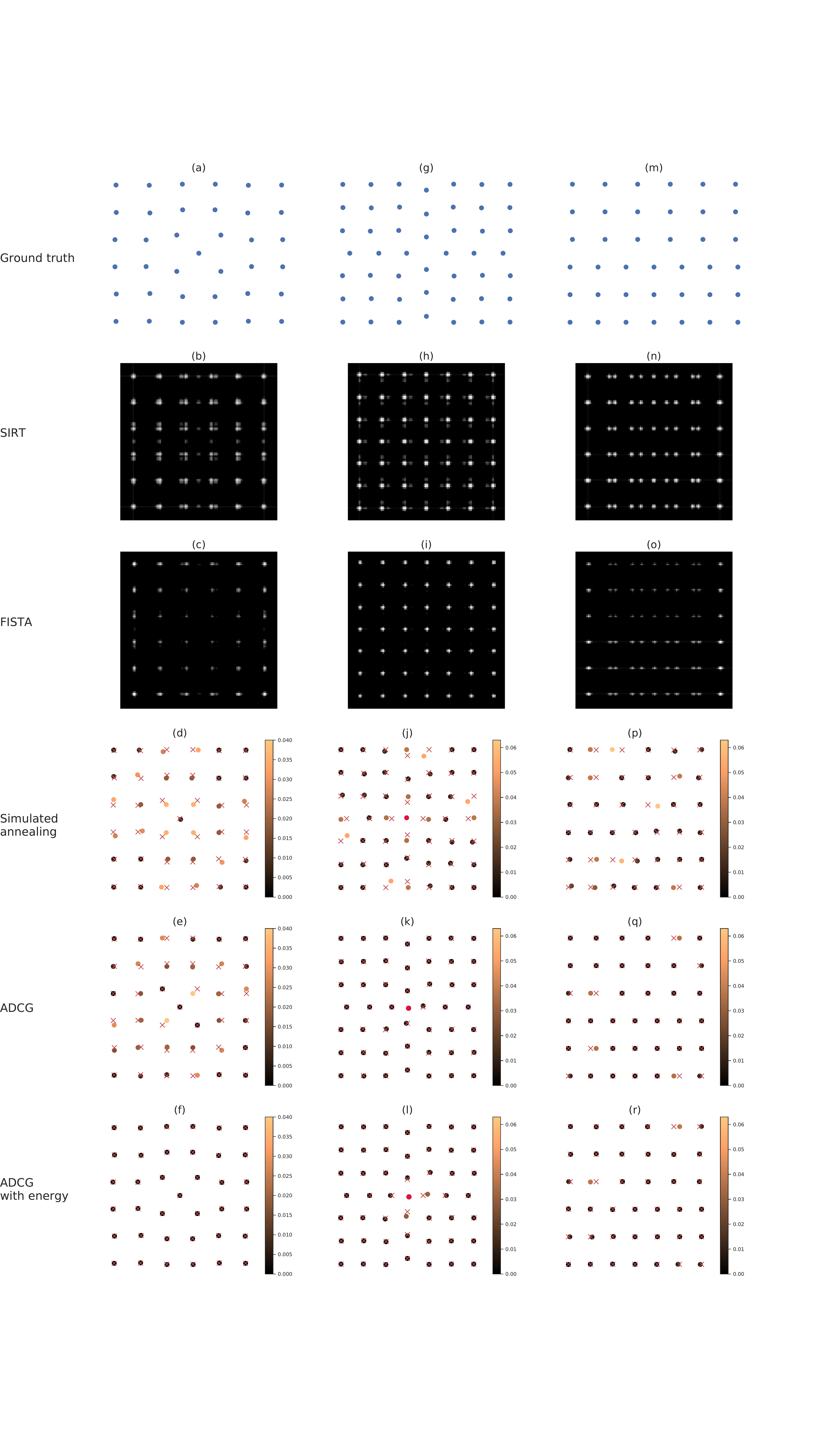}
    \vspace{-3cm}
    \caption{Reconstructions of atomic configurations with (a)--(f) an interstitial point defect, (g)--(l) a vacancy and (m)--(r) an edge dislocation from two projections. For the simulated annealing, ADCG and ADCG with energy reconstructions, atoms are colored according to their Euclidean distance from the ground truth. The ground truth positions are marked with red crosses. In (j)--(l) an extra atom (shown in red) was present in the reconstructions but not in the ground truth.}
    \label{fig:reconstructions}
\end{figure}

Gridding is required for our variant of ADCG (subsection \ref{adcgenergy}, Algorithm \ref{alg:adcg_energy}) at the atom adding step. We found that a coarse discretisation, with less than $1/9^{\text{th}}$ the number of grid points, was already sufficient.

\subsection{Comparison between Reconstructions}
The reconstructions obtained with the different algorithms are shown in Fig.~\ref{fig:reconstructions}. For each reconstruction, data from only two projections were used. Note that two projections is far from sufficient for determining the correct atomic configuration and several different configurations have the same data discrepancy. 


\begin{figure}[t]
    \includegraphics[width=\textwidth]{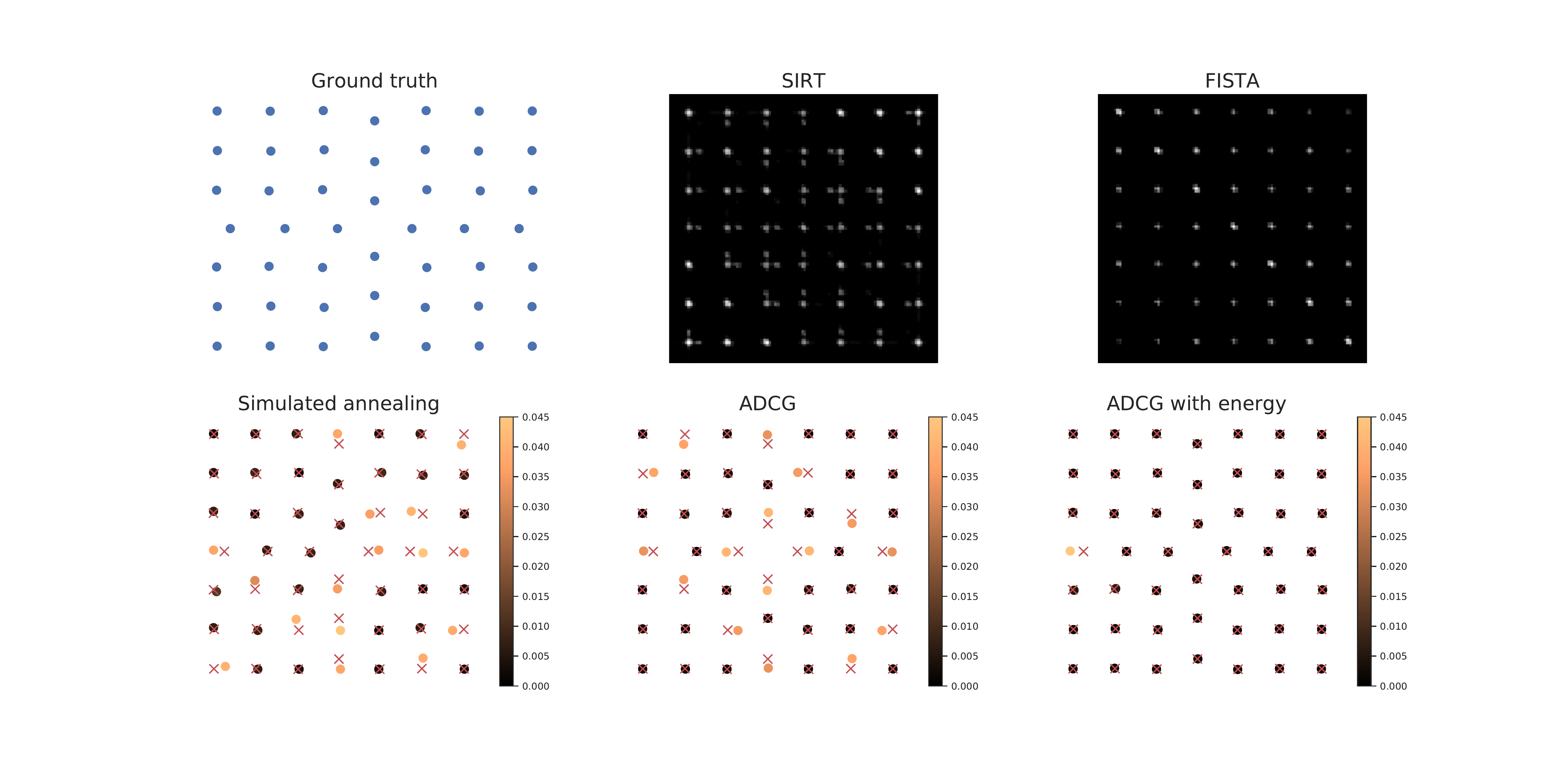}
    \caption{Reconstructions of a vacancy defect from three projections. For the simulated annealing, ADCG and ADCG with energy reconstructions, atoms are colored according to their Euclidean distance from the ground truth. Ground truth positions are marked with red crosses.}
    \label{fig:reconstructions3}
\end{figure}

In the SIRT reconstructions, atom positions were blurred out and none of the defects were resolved. In all cases, the number of intensity peaks was also different from the true number of atoms. Although FISTA reconstructions, which include sparsity constraints on the weights, were less blurry, atoms still occupied more than one pixel. For both these algorithms, additional heuristic post-processing is required to output atom locations. In the edge dislocation case, both algorithms gave rise to a configuration with many more atoms than were present in the ground truth.

\begin{center}
\begin{table}[h!]
\centering
\begin{tabular}{ |p{3cm}||p{1.2cm}|p{1.5cm}||p{1.2cm}|p{1.5cm}||p{1.2cm}|p{1.5cm}||}
 \hline &
 \multicolumn{2}{|c||}{Interstitial defect} & \multicolumn{2}{|c||}{Vacancy (3 projs.~)} & \multicolumn{2}{|c||}{Edge dislocation} \\
 \hline
 & Number of atoms & Mean distance & Number of atoms & Mean distance & Number of atoms & Mean distance \\
 \hline
 Ground truth & 37 & 0.0000 & 48 & 0.0000 & 39 & 0.0000 \\
 \hline
 SIRT & 36 & -- & 49 & -- & 66 & -- \\
 \hline
 FISTA & 36 & -- & 49 & -- & 66 & -- \\
 \hline
 Simulated annealing & 37 & 0.0184 & 48 & 0.0164 & 39 & 0.0159 \\
 \hline
 ADCG & 37 & 0.0138 & 48 & 0.0130 & 39 & 0.0049 \\
 \hline
 ADCG with energy & 37 & 0.0018 & 48 & 0.0024 & 39 & 0.0048\\
 \hline
\end{tabular}
\caption{Number of atoms and mean Euclidean distance from ground truth atoms for reconstructions obtained with different algorithms. Thresholding was used to compute the number of atoms detected in the SIRT and FISTA reconstructions.}
\label{table:metrics}
\end{table}
\end{center}

The discrete simulated annealing algorithm performed better for all configurations. For the interstitial point defect and edge dislocation, the number of atoms in the reconstruction matched that in the ground truth. The positions of most atoms, however, were not resolved correctly. Moreover, the resolution, like in previous algorithms, was limited to the resolution on the detector. We ran the simulated annealing algorithm for comparable times as the ADCG algorithms and picked the solution with the least data discrepancy.

Already the ADCG algorithm for $\alpha = 0$ performed far better than all the previous algorithms. For the configurations with an interstitial point defect and edge dislocation, all but a few atom locations were identified correctly. For the configuration with a vacancy, all atoms were correctly placed. However, an additional atom at the centre of the configuration was placed incorrectly. 

Adding the potential energy (ADCG with energy)  helps to resolve atom positions that were not identified with $\alpha = 0$. For the interstitial point defect and edge dislocation, these reconstructions were the closest to the ground truth. Adding the energy to the configuration with a vacancy moved the atoms near the defect further apart but was not able to correct for the extra atom placed. For this case, we performed an additional experiment with three projections at $0^\circ$, $45^\circ$ and $90^\circ$. These results are shown in Fig.~\ref{fig:reconstructions3}. Taking projections at different angles (e.g.~$0^\circ$, $22.5^\circ$ and $90^\circ$) did not improve results. The defect was still not resolved in the SIRT and FISTA reconstructions. However, the number of atoms in the simulated annealing, ADCG and ADCG with energy reconstructions was correct. Once again, the reconstruction obtained with our algorithm was closer to the ground truth than all other reconstructions, with all but one atom placed correctly. Reconstructions with 3 projections for the interstitial point defect and edge dislocation were not significantly different from those with 2 projections. In Table \ref{table:metrics}, we report the number of atoms detected and (where applicable) the mean Euclidean distance of atoms from the ground truth. Note that for computing the mean distance, we required that the number of atoms detected in the reconstruction matched that in the ground truth. Thresholding with a pre-defined minimum distance between peaks was used to detect atoms in the SIRT and FISTA reconstructions. 

\begin{figure}[t]
  \makebox[\textwidth][c]{\includegraphics[width=1.3\textwidth]{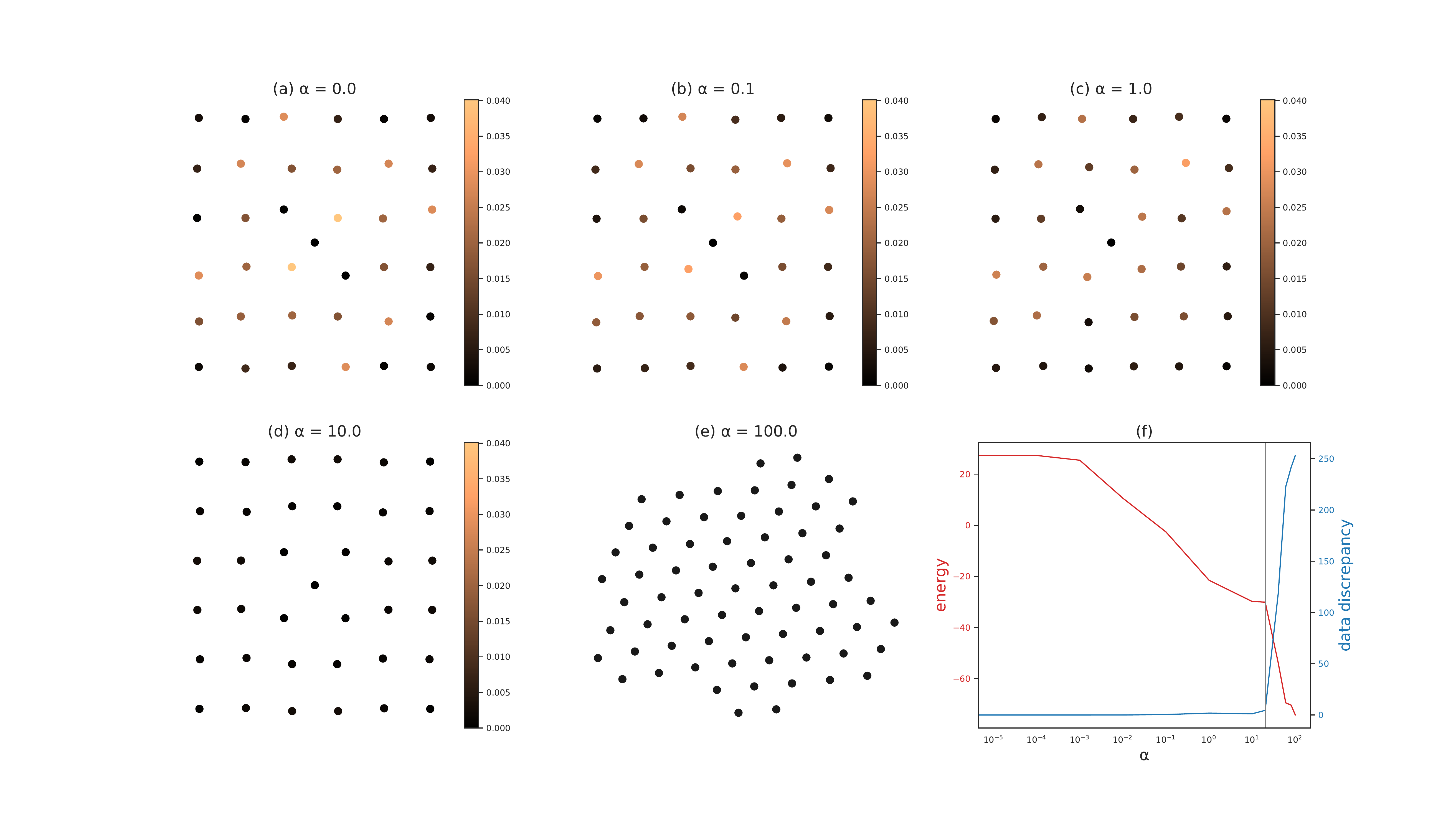}}%
    \caption{(a)-(d): Increasing the weighting of the energy term from $\alpha=0.0$ to $\alpha=10.0$ helps to resolve the correct atomic configuration. The reconstructed atoms are colored according to their Euclidean distance from the atoms in the ground truth. (e) At high values of $\alpha$, the reconstructions have a high data discrepancy and correspond to one of the global minima of the potential energy. (f) From the plots of potential energy and data discrepancy, an optimal value of $\alpha$ (indicated by the grey line) is selected. Increasing $\alpha$ beyond this optimal value leads to a large increase in the data discrepancy due to addition of more atoms.}
    \label{fig:energy}
\end{figure}
\subsection{Effect of Adding Energy to Optimisation}
To resolve atom positions using Algorithm \ref{alg:adcg_energy}, the contribution of the potential energy was increased gradually by increasing $\alpha$. In Fig.~\ref{fig:energy}, we show the effect of adding energy to the optimisation problem. For $\alpha=0$, an initial guess for the true configuration was obtained. This configuration, though data optimal, was not the ground truth. A quantitative measure of this mismatch is the Euclidean distance between the reconstructed atom locations and those in the ground truth. As $\alpha$ was increased, the reconstructions evolved from being data-optimal to being energy-optimal. At a certain value of $\alpha$, the Euclidean distance between reconstructed and ground truth atom locations decreased to zero. Increasing $\alpha$ beyond this point led to a large increase in the data discrepancy term due to the addition of more atoms. For very high values of $\alpha$, the configurations obtained were essentially global minima of the potential energy, such as the honeycomb configuration in Fig.~\ref{fig:energy}(e) for $\alpha=100.0$. An optimal value of the regularisation parameter was selected by increasing $\alpha$ to the point at which more atoms were added to the configuration and a jump in the data discrepancy was observed.

\section{Discussion}\label{discussion}
The results of our numerical experiments demonstrate that algorithms like ADCG, which do not rely on domain discretisation, are better at resolving the defects in the atomic configurations shown in Fig.~\ref{fig:configs}. Moreover, the output from ADCG is a \emph{list of coordinates} and not an image like that of SIRT or FISTA, which requires further post-processing steps to derive the atom locations. Direct access to coordinates can be particularly useful because further analysis, such as strain calculations, often require atom positions as input. 

Adding the potential energy of the atomic configuration to the optimisation problem resulted in reconstructions that were closer to the ground truth.
One challenge of the proposed approach (with or without adding the potential energy) is that the resulting optimisation problem is a non-convex function of the atom locations. The numerical methods we presented are not intentionally designed to escape local minima and are therefore sensitive to their initialisation. To improve this, one important extension would be to also remove atoms from the current configuration, which might make it possible to resolve the vacancy defect in Fig.~\ref{fig:reconstructions} with two projections. More generally, one would need to include suitable features of global optimisation algorithms \cite{pardalos2016advances} that do not compromise ADCG's computational efficiency  (note that we could have adapted simulated annealing to solve \eqref{eq:infoptwithenergy} but using a cooling schedule slow enough to prevent getting trapped in local minima quickly becomes practically infeasible). A related problem is to characterize local and global minimizers of \eqref{eq:infoptwithenergy} to understand which configurations can be uniquely recovered by this approach and which cannot. To process experimental data, it may furthermore be important to analyze the impact of the error caused by the approximate nature of the mathematical models used for data acquisition ($\mathcal{R}, G$) and atomic interaction ($V_{\text{LJ}}$).

\section{Conclusions}\label{conclusions}

In this paper we proposed a novel discrete tomography approach in which the locations of atoms are allowed to vary continuously and their interaction potentials are modeled explicitly. We showed in proof-of-concept numerical studies that such an approach can be better at resolving crystalline defects than image domain discretized or fully discrete algorithms. Furthermore, in situations where atom locations are desired, this approach provides access to the quantity of interest without any additional post-processing. 
For future work, we will extend our numerical studies on this atomic super-resolution approach to larger-sized scenarios in 3D, featuring realistic measurement noise, acquisition geometries, more suitable and accurate physical interaction potentials and different atom types. This will require additional computational effort to scale up our algorithm and will then allow us to work on real electron tomography data of nanocrystals. 

\paragraph{Acknowledgments} This project has received funding from the European Union's Horizon 2020 research and innovation programme under the Marie Sklodowska-Curie grant agreement no.~765604.

\bibliographystyle{splncs04}
\bibliography{mybibliography}

\begin{thebibliography}{10}
\providecommand{\url}[1]{\texttt{#1}}
\providecommand{\urlprefix}{URL }
\providecommand{\doi}[1]{https://doi.org/#1}

\bibitem{van2015astra}
van Aarle, W., Palenstijn, W.J., De~Beenhouwer, J., Altantzis, T., Bals, S.,
  Batenburg, K.J., Sijbers, J.: The astra toolbox: A platform for advanced
  algorithm development in electron tomography. Ultramicroscopy  \textbf{157},
  35--47 (2015)

\bibitem{adler2017operator}
Adler, J., Kohr, H., Oktem, O.: Operator discretization library (odl). Software
  available from https://github. com/odlgroup/odl  (2017)

\bibitem{alberti2019dynamic}
Alberti, G.S., Ammari, H., Romero, F., Wintz, T.: Dynamic spike superresolution
  and applications to ultrafast ultrasound imaging. SIAM Journal on Imaging
  Sciences  \textbf{12}(3),  1501--1527 (2019)

\bibitem{alpers2006}
Alpers, A., Gritzmann, P.: On stability, error correction, and noise
  compensation in discrete tomography. SIAM Journal on Discrete Mathematics
  \textbf{20}(1),  227--239 (2006)

\bibitem{anderson2008general}
Anderson, J.A., Lorenz, C.D., Travesset, A.: General purpose molecular dynamics
  simulations fully implemented on graphics processing units. Journal of
  computational physics  \textbf{227}(10),  5342--5359 (2008)

\bibitem{baake2006}
Baake, M., Huck, C., Gritzmann, P., Langfeld, B., Lord, K.: Discrete tomography
  of planar model sets. Acta Crystallographica A  \textbf{62}(6),  419--433
  (2006)

\bibitem{batenburg2006}
Batenburg, K.J.: A network flow algorithm for reconstructing binary images from
  discrete x-rays. Journal of Mathematical Imaging and Vision  \textbf{27}(2),
  175--191 (2006)

\bibitem{batenburg2009}
Batenburg, K.J., Sijbers, J.: Generic iterative subset algorithms for discrete
  tomography. Discrete Applied Mathematics  \textbf{157}(3),  438--451 (2009)

\bibitem{beck2009fast}
Beck, A., Teboulle, M.: A fast iterative shrinkage-thresholding algorithm for
  linear inverse problems. SIAM journal on imaging sciences  \textbf{2}(1),
  183--202 (2009)

\bibitem{bitzek2006structural}
Bitzek, E., Koskinen, P., G{\"a}hler, F., Moseler, M., Gumbsch, P.: Structural
  relaxation made simple. Physical review letters  \textbf{97}(17),  170201
  (2006)

\bibitem{boyd2017alternating}
Boyd, N., Schiebinger, G., Recht, B.: The alternating descent conditional
  gradient method for sparse inverse problems. SIAM Journal on Optimization
  \textbf{27}(2),  616--639 (2017)

\bibitem{bredies2013inverse}
Bredies, K., Pikkarainen, H.K.: Inverse problems in spaces of measures. ESAIM:
  Control, Optimisation and Calculus of Variations  \textbf{19}(1),  190--218
  (2013)

\bibitem{candes2014towards}
Cand{\`e}s, E.J., Fernandez-Granda, C.: Towards a mathematical theory of
  super-resolution. Communications on pure and applied Mathematics
  \textbf{67}(6),  906--956 (2014)

\bibitem{chen2013three}
Chen, C.C., Zhu, C., White, E.R., Chiu, C.Y., Scott, M., Regan, B., Marks,
  L.D., Huang, Y., Miao, J.: Three-dimensional imaging of dislocations in a
  nanoparticle at atomic resolution. Nature  \textbf{496}(7443),  74--77 (2013)

\bibitem{denoyelle2019sliding}
Denoyelle, Q., Duval, V., Peyr{\'e}, G., Soubies, E.: The sliding frank-wolfe
  algorithm and its application to super-resolution microscopy. Inverse
  Problems  (2019)

\bibitem{foucart2013mathematical}
Foucart, S., Rauhut, H.: A mathematical introduction to compressive sensing.
  Applied and Numerical Harmonic Analysis, Birkh\"{a}user/Springer, New York
  (2013)

\bibitem{frank1956algorithm}
Frank, M., Wolfe, P.: An algorithm for quadratic programming. Naval research
  logistics quarterly  \textbf{3}(1-2),  95--110 (1956)

\bibitem{frenkel2001understanding}
Frenkel, D., Smit, B.: Understanding molecular simulation: from algorithms to
  applications, vol.~1. Elsevier (2001)

\bibitem{gardner1997}
Gardner, R.J., Gritzmann, P.: Discrete tomography: Determination of finite sets
  by x-rays. Transactions of the American Mathematical Society
  \textbf{349}(6),  2271--2295 (1997)

\bibitem{glaser2015strong}
Glaser, J., Nguyen, T.D., Anderson, J.A., Lui, P., Spiga, F., Millan, J.A.,
  Morse, D.C., Glotzer, S.C.: Strong scaling of general-purpose molecular
  dynamics simulations on gpus. Computer Physics Communications  \textbf{192},
  97--107 (2015)

\bibitem{goris2015measuring}
Goris, B., De~Beenhouwer, J., De~Backer, A., Zanaga, D., Batenburg, K.J.,
  Sannchez-Iglesias, A., Liz-Marzan, L.M., Van~Aert, S., Bals, S., Sijbers, J.,
  et~al.: Measuring lattice strain in three dimensions through electron
  microscopy. Nano letters  \textbf{15}(10),  6996--7001 (2015)

\bibitem{jaggi2013revisiting}
Jaggi, M.: Revisiting frank-wolfe: Projection-free sparse convex optimization.
  In: ICML (1). pp. 427--435 (2013)

\bibitem{kak2002principles}
Kak, A.C., Slaney, M., Wang, G.: Principles of computerized tomographic
  imaging. Medical Physics  \textbf{29}(1),  107--107 (2002)

\bibitem{montroll1956theory}
Montroll, E.W.: Theory of the vibration of simple cubic lattices with nearest
  neighbor interactions. In: Proceedings of the Third Berkeley Symposium on
  Mathematical Statistics and Probability. vol.~3, pp. 209--246. Univ of
  California Press (1956)

\bibitem{nelder1965simplex}
Nelder, J.A., Mead, R.: A simplex method for function minimization. The
  computer journal  \textbf{7}(4),  308--313 (1965)

\bibitem{pardalos2016advances}
Pardalos, P.M., Zhigljavsky, A., {\v{Z}}ilinskas, J.: Advances in stochastic
  and deterministic global optimization. Springer (2016)

\bibitem{parikh2014proximal}
Parikh, N., Boyd, S., et~al.: Proximal algorithms. Foundations and
  Trends{\textregistered} in Optimization  \textbf{1}(3),  127--239 (2014)

\bibitem{poon2019multidimensional}
Poon, C., Peyr{\'e}, G.: Multidimensional sparse super-resolution. SIAM Journal
  on Mathematical Analysis  \textbf{51}(1),  1--44 (2019)

\bibitem{rez2013three}
Rez, P., Treacy, M.M.: Three-dimensional imaging of dislocations. Nature
  \textbf{503}(7476),  E1--E1 (2013)

\bibitem{van2011three}
Van~Aert, S., Batenburg, K.J., Rossell, M.D., Erni, R., Van~Tendeloo, G.:
  Three-dimensional atomic imaging of crystalline nanoparticles. Nature
  \textbf{470}(7334), ~374 (2011)

\bibitem{van2012big}
Van~Dyck, D., Jinschek, J.R., Chen, F.R.: ‘big bang’tomography as a new
  route to atomic-resolution electron tomography. Nature  \textbf{486}(7402),
  243--246 (2012)

\bibitem{virtanen2019scipy}
Virtanen, P., Gommers, R., Oliphant, T.E., Haberland, M., Reddy, T.,
  Cournapeau, D., Burovski, E., Peterson, P., Weckesser, W., Bright, J.,
  et~al.: Scipy 1.0--fundamental algorithms for scientific computing in python.
  arXiv preprint arXiv:1907.10121  (2019)

\end{thebibliography}

\end{document}